\documentclass[12pt]{article}
\usepackage{amsmath}
\usepackage{amsfonts}
\usepackage{amssymb}
\usepackage{mathrsfs}
\usepackage{graphicx}
\usepackage{epsfig}
\usepackage{amsthm}
\usepackage{latexsym}
\usepackage{cite}

\newtheorem{theorem}{Theorem}[section]
\newtheorem{defi}[theorem]{Definition}
\newtheorem{lemma}[theorem]{Lemma}
\newtheorem{remark}[theorem]{Remark}

\def\proof {{\noindent\bf Proof.}\quad}
\def \endproof{\hfill$\Box$\vspace {3mm}}

\begin{document}

\setlength{\oddsidemargin}{0cm} \setlength{\evensidemargin}{0cm}
\baselineskip=20pt

\title{\textbf{Brake subharmonic solutions of first order Hamiltonian systems}\footnote{Partially
supported by the NNSFC (10531050,10621101) and 973 Program of STM
(2006CB805903)}}
\author{ Chong Li\thanks{E-mail: plumechong@yahoo.com.cn}\quad\quad Chungen Liu
\thanks{E-mail: liucg@nankai.edu.cn}
\\School of Mathematical Sciences and LPMC, \\Nankai University, Tianjin
300071, P.R. China.}
\date{}

\maketitle

\begin{center}
 \baselineskip=18pt
 \begin{minipage}{5in}
{\textbf {Abstract}}\quad  In this paper, we mainly use the Galerkin
approximation method and the iteration inequalities of the
$L$-Maslov type index theory in \cite{liu2, liu5} to study the
properties of brake subharmonic solutions for the first order
non-autonomous Hamiltonian systems. We prove that when the positive
integers $j$ and $k$ satisfies the certain conditions, there exists
a $jT$-periodic nonconstant brake solution $z_{j}$ such that $z_{j}$
and $z_{kj}$ are distinct.
\medskip


\medskip

{\textbf{2000 MR Subject Classification}}\quad 58F05, 58E05, 34C25,
58F10
\end{minipage}
\end{center}
\section{\large\textbf{Introduction and the Main Results}}
\renewcommand{\theequation}{\arabic{section}.\arabic{equation}}

\qquad In this paper, we consider the first order non-autonomous
Hamiltonian systems
\begin{eqnarray}\label{e1.1}
\dot{z}(t)=J\nabla H(t,z(t)),\; \forall z \in \mathbb{R}^{2n},\;
\forall t \in \mathbb{R},
\end{eqnarray}
where  $J=\left(\begin{matrix}0 & -I_{n}
\\ I_{n} & 0 \end{matrix} \right)$ is the standard symplectic matrix,
  $I_{n}$ is the
unit matrix of order $n$, $H\in C^{2}(\mathbb{R} \times
{\mathbb{R}}^{2n}, \mathbb{R})$ and $\nabla H(t,z)$ is the gradient
of $H(t,z)$ with respect to the space variable $z$. We denote the
standard norm and inner product in $\mathbb{R}^{2n}$ by $|\cdot|$
and $(\cdot,\cdot)$, respectively.

Suppose that $H(t,z)=\frac{1}{2} (\hat{B}(t)z,z)+\hat{H}(t,z)$ and
$H\in C^{2}(\mathbb{R} \times {\mathbb{R}}^{2n}, \mathbb{R})$
satisfies the following conditions:

(H1) $\hat{H}(T+t, z)=\hat{H}(t, z)$, for all $z\in
\mathbb{R}^{2n}$, $t\in \mathbb{R}$,

(H2) $\hat{H}(t, z)=\hat{H}(-t, Nz)$, for all $z\in
\mathbb{R}^{2n}$, $t\in \mathbb{R}$, $N=\left(\begin{matrix} -I_n &
0\\0 & I_n\end{matrix}\right)$,

(H3) $\hat{H}''(t, z)>0$, for all $z\in \mathbb{R}^{2n}\backslash
\{0\}$, $t\in \mathbb{R}$,

(H4) $\hat{H}(t,z)\geq 0$, for all $z\in \mathbb{R}^{2n}$, $t\in
\mathbb{R}$,

(H5) $\hat{H}(t,z)=o (|z|^{2})$ at $z=0$,

(H6) There is a $\theta\in (0,1/2)$ and $\bar{r}> 0$ such that
\begin{eqnarray*}
0< \frac{1}{\theta}\hat{H}(t,z)\leq (z,\nabla \hat{H}(t,z)), \mbox{
for all } z\in \mathbb{R}^{2n}, |z|\geq \bar{r}, t\in \mathbb{R},
\end{eqnarray*}

(H7) $\hat{B}(t)$ is a symmetrical continuous matrix,
$|\hat{B}|_{C^{0}}\leq \beta_0$ for some $\beta_0 >0$, and
$\hat{B}(t)$ is a semi-positively definite for all $t\in
\mathbb{R}$,

(H8) $\hat{B}(T+t)=\hat{B}(t)=\hat{B}(-t)$,
$\hat{B}(t)N=N\hat{B}(t)$, for all $t\in \mathbb{R}$.

Recall that a $T$-periodic solution $(z, T)$ of (\ref{e1.1}) is
called {\it brake solution} if $z(t+T)=z(t)$ and $z(t)=Nz(-t)$, the
later is equivalent to $z(T/2+t)=Nz(T/2-t)$, in this time $T$ is
called the {\it brake period} of $z$. Up to the authors' knowledge,
H. Seifert firstly studied brake orbits in second order autonomous
Hamiltonian systems in \cite{se} of 1948. Since then many studies
have been carried out for brake orbits of first order and second
order Hamiltonian systems. For the minimal periodic problem,
multiple existence results about brake orbits for the Hamiltonian
systems and more details on brake orbits one can refer the papers
\cite{ABL1, Be1, BG, Bol, BolZ, Gro, GZ1, Ha1, liu5, lo2, Ra1, Sz}
and the references therein. S. Bolotin proved first in \cite{Bol}
(also see \cite{BolZ}) of 1978 the existence of brake orbits in
general setting. K. Hayashi in \cite{Ha1}, H. Gluck and W. Ziller in
\cite{GZ1}, and V. Benci in \cite{Be1} in 1983-1984 proved the
existence of brake orbits of second order Hamiltonian systems under
certain conditions. In 1987, P. Rabinowitz in \cite{Ra1} proved the
existence of brake orbits of first order Hamiltonian systems. In
1987, V. Benci and F. Giannoni gave a different proof of the
existence of one brake orbit in \cite{BG}. In 1989, A. Szulkin in
\cite{Sz} proved the existence of brake orbits of first order
Hamiltonian systems under the $\sqrt{2}$-pinched condition. E. van
Groesen in \cite{Gro} of 1985 and A. Ambrosetti, V. Benci, Y. Long
in \cite{ABL1} of 1993 also proved the multiplicity result about
brake orbits for the second order Hamiltonian systems under
different pinching conditions. Without pinching conditions, in
\cite{lo2} (2006) Y. Long, D. Zhang and C. Zhu proved that there
exist at least two geometrically distinct brake orbits in every
bounded convex symmetric domain in $\mathbb{R}^{n}$ for $n\geq 2$.
Recently, C. Liu and D. Zhang in \cite{liu5} proved that there exist
at least $[n/2]+1$ geometrically distinct brake orbits in every
bounded convex symmetric domain in $\mathbb{R}^{n}$ for $n\geq 2$,
and there exist at least $n$ geometrically distinct brake orbits on
nondegenerate domain.

For the non-autonomous Hamiltonian systems, for periodic boundary
(brake solution) problems, since the Hamiltonian function $H$ is
$T$-periodic in the time variable $t$, if the system (\ref{e1.1})
has a $T$-periodic solution $(z_{1}, T)$, one hopes to find the
$jT$-periodic solution $(z_{j}, jT)$ for integer $j\geq 1$, for
example, $(z_1, jT)$ itself is $jT$-periodic solution. The
subharmonic solution problem asks when the solutions $z_{1}$ and
$z_{j}$ are distinct. More precisely, in the case of brake
solutions, $z_{1}$ and $z_{j}$ are {\it distinct} if
$\frac{kT}{2}*z_{1}(\cdot)\equiv z_{1}(\frac{kT}{2}+\cdot)\neq
z_{j}(\cdot)$ for any integer $k$. In other word, $z_j(t)\ne
z_1(t))$ and $z_j(t)\ne z_1(T/2+t)$ for $t\in [0,T]$. In below we
remind that the $L_{0}$-indices of the two solutions $z_{1}$ and
$(kT)*z_{1}$ for any $k\in \mathbb{Z}$ in the interval $[0, T/2]$
are the same. In this paper, we first consider the brake subharmonic
solution problem. We state the main results of this paper.

\begin{theorem}\label{t1.1}
Suppose that $H\in
C^{2}(\mathbb{R} \times {\mathbb{R}}^{2n}, \mathbb{R})$ satisfies
(H1)-(H8), then for each integer $1\leq j< 2\pi/\beta_0 T$, there is
a $jT$-periodic nonconstant brake solution $z_{j}$ of (\ref{e1.1})
such that $z_{j}$ and $z_{kj}$ are distinct for $k\geq 5$ and $kj<
2\pi/\beta_0 T$. Furthermore, $\{z_{k^{p}}| p\in \mathbb{N}\}$ is a
pairwise distinct brake solution sequence of (\ref{e1.1}) for $k\geq
5$ and $1\leq k^{p}< 2\pi/\beta_0 T$.
\end{theorem}

Especially, if $\hat{B}(t)\equiv 0$, then $2\pi/ \beta_0 T=+\infty$.
Therefore, one can state the following theorem.

\begin{theorem}\label{t1.2}
Suppose that $H\in C^{2}(\mathbb{R} \times {\mathbb{R}}^{2n},
\mathbb{R})$ with $\hat{B}(t)\equiv 0$ satisfies (H1)-(H6), then for
each integer $j\geq 1$, there is a $jT$-periodic nonconstant brake
solution $z_{j}$ of (\ref{e1.1}). Furthermore, given any integers
$j\geq 1$ and $k\geq 5$, $z_{j}$ and $z_{kj}$ are distinct brake
solutions of (\ref{e1.1}), in particularly, $\{z_{k^{p}}| p\in
\mathbb{N}\}$ is a pairwise distinct brake solution sequence of
(\ref{e1.1}).
\end{theorem}

The first result on subharmonic periodic solutions for the
Hamiltonian systems $\dot{z}(t)=J \nabla H(t, z(t))$, where $z\in
\mathbb{R}^{2n}$ and $H(t, z)$ is $T$-periodic in $t$, was obtained
by P. Rabinowitz in his pioneer work \cite{ra2}. Since then, many
new contributions have appeared. See for example \cite{ek, ekho,
liu4, lo1, si} and the references therein. Especially, in
\cite{ekho}, I. Ekeland and H. Hofer proved that under a strict
convex condition and a superquadratic condition, the Hamiltonian
system $\dot{z}(t)=J \nabla H(t, z(t))$ possesses subharmonic
solution $z_{k}$ for each integer $k \geq 1$ and all of these
solutions are pairwise geometrically distinct. In \cite{liu4}, the
second author of this paper obtained a result of subharmonic
solutions for the non-convex case by using the Maslov-type index
iteration theory. We notice that in \cite{an} T. An wants to improve
the result of \cite{liu4}, but there is a gap in his proof when
applying Theorem 2.6 there to prove his Theorem 1.3. Precisely, the
formula (2.17) in \cite{an} is not true since the middle term
$i_T(\hat B)+\nu_T(\hat B)+1$ should be $i_{kT}(\hat
B)+\nu_{kT}(\hat B)+1$. Up to the authors' knowledge, Theorem
\ref{t1.1} and \ref{t1.2} are the first results for the brake
subharmonic solution problem for the time being.

The main ingredient in proving Theorem \ref{t1.1} and \ref{t1.2} is
to transform the brake solution problem into the $L_{0}$-boundary
problem:
\begin{eqnarray}\label{e1.2}
\left\{\begin{array}{ll}&\dot{z}(t)=J\nabla H(t,z(t)),\; \forall z
\in \mathbb{R}^{2n},\; \forall t \in [0, T/2], \\
&z(0)\in L_{0},\; z(T/2)\in L_{0},
 \end{array}\right.
 \end{eqnarray}
where $L_{0}=\{0\} \oplus \mathbb{R}^{n} \in \Lambda (n)$.
$\Lambda(n)$ is the set of all linear Lagrangian subspaces in
$(\mathbb{R}^{2n},\omega_{0})$, here the standard symplectic form is
defined by $\omega_{0}=\sum\limits_{i=1}^{n}dx_{i}\wedge dy_{i}$. A
Lagrangian subspace $L$ of ${\mathbb{R}}^{2n}$ is an $n$ dimensional
subspace satisfying $\omega_{0}|_{L}=0$.

\begin{lemma}\label{l1.1}
Suppose the Hamiltonian function $H$
satisfying conditions (H1), (H2) and (H8). If $(z, T/2)$ is a
solution of the problem (\ref{e1.2}), then $(\tilde{z}, T)$ is a
$T$-periodic solution of the Hamiltonian system (1.1) satisfying the
brake condition $\tilde{z}(T/2+t)=N\tilde{z}(T/2-t)$, where
$\tilde{z}$ is defined by
\begin{eqnarray*}
\tilde{z}(t)=\left\{\begin{array}{ll}z(t),\;& t\in
[0,T/2],\\Nz(T-t),\;&t\in (T/2,T].\end{array}\right.
\end{eqnarray*}
\end{lemma}

\proof It is easy to see that $\tilde{z}(t)$ is continuous in the
interval $[0, T]$. By direct computation,
\begin{eqnarray*}
\dot{\tilde{z}}(t+T/2)=-N\dot{z}(T/2-t) =JN\nabla
H(T/2-t,z(T/2-t))\\
=J\nabla H(t+T/2,Nz(T/2-t))=J\nabla H(t+T/2,\tilde{z}(t+T/2)).
\end{eqnarray*}
So $(\tilde{z}, T)$ is a $T$-periodic solution of the Hamiltonian
system (\ref{e1.1}). The brake condition is satisfied by the
definition of $\tilde{z}$. The proof of Lemma \ref{l1.1} is
complete.\endproof

By this observation, we then use the Galerkin approximation methods
to get a critical point of the action functional which is also a
solution of (\ref{e3.1}) with a suitable $L_{0}$-index estimate, see
Theorem \ref{t3.1} below. The $L$-Maslov type index theory for any
$L\in \Lambda(n)$ was studied in \cite{liu3} by the algebraic
methods. In \cite{lo2}, Y. Long, D.Zhang and C. Zhu established two
indices $\mu_1(\gamma)$ and $\mu_2(\gamma)$ for the fundamental
solution $\gamma$ of a linear Hamiltonian system by the methods of
functional analysis which are special cases of the $L$-Maslov type
index $i_L(\gamma)$ for Lagrangian subspaces $L_0=\{0\}\oplus
{\mathbb R}^n$ and $L_1={\mathbb R}^n\oplus \{0\}$ up to a constant
$n$. In order to prove Theorem \ref{t1.1} and \ref{t1.2}, we need to
consider the problem (\ref{e3.1}). The iteration theory of the
$L_{0}$-Maslov type index theory developed in \cite{liu2} and
\cite{liu5}, then help us to distinguish solutions $z_{j}$ from
$z_{kj}$ in Theorem \ref{t1.1} and \ref{t1.2}.

This paper is divided into 3 sections. In section 2, we give a brief
introduction to the Maslov-type index theory for symplectic paths
with Lagrangian boundary conditions and an iteration theory for the
$L_{0}$-Maslov type index theory. In section 3, we give a proof of
Theorem \ref{t1.1} and \ref{t1.2}.

\section{\large\textbf{Preliminaries}}
\setcounter{equation}{0}
\renewcommand{\theequation}{\arabic{section}.\arabic{equation}}

\qquad In this section, we briefly recall the Maslov-type index
theory for symplectic paths with Lagrangian boundary conditions and
an iteration theory for the $L_{0}$-Maslov type index theory. All
the details can be found in \cite{liu1, liu2, liu3, liu5}.

We denote the $2n$-dimensional symplectic group $Sp(2n)$ by
\begin{eqnarray*}
Sp(2n)=\{M\in \mathscr{L}(\mathbb{R}^{2n})| M^{T}JM=J\},
\end{eqnarray*}
where $\mathscr{L}(\mathbb{R}^{2n})$ is the set of all real
$2n\times 2n$ matrices, $M^{T}$ is the transpose of matrix $M$.
Denote by $\mathscr{L}_{s}(\mathbb{R}^{2n})$ the subset of
$\mathscr{L}(\mathbb{R}^{2n})$ consisting of symmetric matrices. And
denote the symplectic path space by
\begin{eqnarray*}
\mathcal {P}(2n)=\{\gamma \in C([0, 1], Sp(2n)) |
\gamma(0)=I_{2n}\}.
\end{eqnarray*}
We write a symplectic path $\gamma \in \mathcal {P}(2n)$ in the
following form
\begin{eqnarray}\label{e2.1}
\gamma(t)=\left(\begin{matrix}S(t) & V(t)
\\ T(t) & U(t) \end{matrix} \right),
\end{eqnarray}
 where $S(t)$, $T(t)$, $V(t)$, $U(t)$ are $n\times n$ matrices. The
 $n$ vectors come from the column of the matrix $\left(\begin{matrix}V(t)
 \\ U(t) \end{matrix} \right)$ are linear
 independent and they span a Lagrangian subspace of $(\mathbb{R}^{2n},
 \omega_{0})$. Particularly, at $t=0$, this Lagrangian subspace is
 $L_{0}=\{0\}\oplus \mathbb{R}^{2n}$.

\begin{defi} {\rm(see \cite{liu3})} \label{d2.1}
We define the $L_{0}$-nullity of any symplectic path $\gamma \in
\mathcal {P}(2n)$ by
\begin{eqnarray*}
\nu_{L_{0}}(\gamma)\equiv\dim\ker_{L_{0}}(\gamma(1)):=\dim\ker
V(1)=n- \mbox{{\rm rank }}V(1)
\end{eqnarray*}
with the $n\times n$ matrix function $V(t)$ defined in (\ref{e2.1}).
\end{defi}

For $L_{0}=\{0\}\oplus \mathbb{R}^{n}$, We define the following
subspaces of $Sp(2n)$ by
\begin{eqnarray*}
Sp(2n)^{*}_{L_{0}}&=&\{M \in Sp(2n) | \det V_{M}\neq 0\},\\
Sp(2n)^{0}_{L_{0}}&=&\{M \in Sp(2n) | \det V_{M}=0\},\\
Sp(2n)^{\pm}_{L_{0}}&=&\{M \in Sp(2n) | \pm\det V_{M}>0\},
\end{eqnarray*}
where $M=\left(\begin{matrix}S_{M} & V_{M}
\\ T_{M} & U_{M} \end{matrix} \right)$ and
$Sp(2n)_{L_{0}}^{*}=Sp(2n)_{L_{0}}^{+}\cup Sp(2n)_{L_{0}}^{-}$. And
denote two subsets of $\mathcal {P}(2n)$ by
\begin{eqnarray*}
\mathcal {P}(2n)^{*}_{L_{0}}&=&\{\gamma \in \mathcal {P}(2n) | \nu_{L_{0}}(\gamma)=0\},\\
\mathcal {P}(2n)^{0}_{L_{0}}&=&\{\gamma \in \mathcal {P}(2n) |
\nu_{L_{0}}(\gamma)>0\}.
\end{eqnarray*}
\qquad We note that $\mbox{rank }\left(\begin{matrix}V(t)
 \\ U(t) \end{matrix} \right)=n$, so the complex matrix $U(t)\pm
 \sqrt{-1}V(t)$ is invertible. We define a complex matrix function
 by
 \begin{eqnarray*}
 Q(t)=(U(t)-\sqrt{-1}V(t))(U(t)+\sqrt{-1}V(t))^{-1}.
 \end{eqnarray*}
 It is easy to see that the matrix $Q(t)$ is a unitary matrix for
 any $t\in [0,1]$. We define
 \begin{eqnarray*}
 M_{+}=\left(\begin{matrix}0 & I_{n}
\\ -I_{n} & 0 \end{matrix} \right),\;
M_{-}=\left(\begin{matrix}0 & J_{n}
\\ -J_{n} & 0 \end{matrix} \right),\;
J_{n}=\mbox{diag }(-1, 1, \cdots, 1).
\end{eqnarray*}

For a path $\gamma \in \mathcal {P}(2n)^{*}_{L_{0}}$, we first
adjoin it with a simple symplectic path starting from $J=-M_{+}$,
that is, we define a symplectic path by
\begin{eqnarray*}
\tilde{\gamma}(t)=\left\{\begin{array}{ll}I\cos
\frac{(1-2t)\pi}{2}+J\sin \frac{(1-2t)\pi}{2},\; & t \in [0, 1/2], \\
\gamma (2t-1),\; & t\in [1/2, 1].
 \end{array}\right.
 \end{eqnarray*}
 Then we choose a symplectic path $\beta (t)$ in $Sp(2n)^{*}_{L_{0}}$
 starting from $\gamma (1)$ and ending at $M_{+}$ or $M_{-}$
 according to $\gamma(1)\in Sp(2n)_{L_{0}}^{+}$ or $\gamma(1)\in
 Sp(2n)_{L_{0}}^{-}$, respectively. We now define a joint path by
\begin{eqnarray*}
\bar{\gamma}(t)=\beta \ast
\tilde{\gamma}:=\left\{\begin{array}{ll}\tilde{\gamma}(2t),\; & t \in [0, 1/2], \\
\beta (2t-1),\; & t\in [1/2, 1].
 \end{array}\right.
 \end{eqnarray*}
 By the definition, we see that the symplectic path $\bar{\gamma}$
 starting from $-M_{+}$ and ending at either $M_{+}$ or $M_{-}$. As
 above, we define
 \begin{eqnarray*}
 \bar{Q}(t)=(\bar{U}(t)-\sqrt{-1}\bar{V}(t))(\bar{U}(t)+\sqrt{-1}\bar{V}(t))^{-1},
 \end{eqnarray*}
 for $\bar{\gamma}(t)=\left(\begin{matrix}\bar{S}(t) & \bar{V}(t)
\\ \bar{T}(t) & \bar{U}(t) \end{matrix} \right)$. We can choose a
continuous function $\bar{\Delta}(t)$ in [0,1] such that
\begin{eqnarray*}
\mbox{det }\bar{Q}(t)=e^{2\sqrt{-1}\bar{\Delta}(t)}.
\end{eqnarray*}
By the above arguments, we see that the number
$\frac{1}{\pi}(\bar{\Delta}(1)-\bar{\Delta}(0))\in \mathbb{Z}$ and
it does not depend on the choice of the  function $\bar{\Delta}(t)$.

\begin{defi} {\rm(see \cite{liu3})} \label{d2.2}
For a symplectic path
$\gamma \in \mathcal {P}(2n)^{*}_{L_{0}}$, we define the
$L_{0}$-index of $\gamma$ by
\begin{eqnarray*}
i_{L_{0}}(\gamma)=\frac{1}{\pi}(\bar{\Delta}(1)-\bar{\Delta}(0)).
\end{eqnarray*}
\end{defi}

\begin{defi} {\rm(see \cite{liu3})} \label{d2.3}
For a symplectic path
$\gamma \in \mathcal {P}(2n)^{0}_{L_{0}}$, we define the
$L_{0}$-index of $\gamma$ by
\begin{eqnarray*}
i_{L_{0}}(\gamma)=\inf \{i_{L_{0}}(\tilde{\gamma}) |
\tilde{\gamma}\in \mathcal {P}(2n)^{*}_{L_{0}}, \mbox{and
}\tilde{\gamma} \mbox{ is sufficiently close to }\gamma\}.
\end{eqnarray*}
\end{defi}

We know that $\Lambda(n)=U(n)/O(n)$, this means that for any linear
subspace $L\in \Lambda(n)$, there is an orthogonal symplectic matrix
$P=\left(\begin{matrix}A & -B \\ B & A \end{matrix} \right)$ with
$A\pm \sqrt{-1}B \in U(n)$, the unitary matrix, such that
$PL_{0}=L$. $P$ is uniquely determined by $L$ up to an orthogonal
matrix $C \in O(n)$. It means that for any other choice $P'$
satisfying above conditions, there exists
a matrix $C \in O(n)$ such that $P'=P\left(\begin{matrix}C & 0 \\
0 & C
\end{matrix} \right)$ (see \cite{mc}). We define the conjugated
symplectic path $\gamma_{c}\in \mathcal {P}(2n)$ of $\gamma$ by
$\gamma_{c}(t)=P^{-1}\gamma(t)P$.

\begin{defi} {\rm(see \cite{liu3})} \label{d2.4}
We define the $L$-nullity
of any symplectic path $\gamma \in \mathcal {P}(2n)$ by
\begin{eqnarray*}
\nu_{L}(\gamma)\equiv\dim\ker_{L}(\gamma(1)):=\dim\ker
V_{c}(1)=n-\mbox{{\rm rank }}V_{c}(1),
\end{eqnarray*}
where the $n\times n$ matrix function $V_{c}(t)$ is defined in
(\ref{e2.1}) with the symplectic path $\gamma$ replaced by
$\gamma_{c}$, i.e., $\gamma_{c}(t)=\left(\begin{matrix}S_{c}(t) &
V_{c}(t)
\\ T_{c}(t) & U_{c}(t) \end{matrix} \right)$.
\end{defi}

\begin{defi} {\rm(see \cite{liu3})} \label{d2.5}
For a symplectic path
$\gamma \in \mathcal {P}(2n)$, we define the $L$-index of $\gamma$
by
\begin{eqnarray*}
i_{L}(\gamma)=i_{L_{0}}(\gamma_{c}).
\end{eqnarray*}
\end{defi}

In the case of linear Hamiltonian systems
\begin{eqnarray}\label{e2.2}
\dot{y}=JB(t)y,\quad \forall y\in \mathbb{R}^{2n},
\end{eqnarray}
where $B\in C(\mathbb{R}, \mathscr{L}_{s}(\mathbb{R}^{2n}))$. Its
fundamental solution $\gamma=\gamma_{B}$ is a symplectic path
starting from identity matrix $I_{2n}$, i.e., $\gamma=\gamma_{B} \in
\mathcal {P}(2n)$. We denote by
\begin{eqnarray*}
i_{L}(B)=i_{L}(\gamma_{B}),\; \nu_{L}(B)=\nu_{L}(\gamma_{B}).
\end{eqnarray*}

\begin{theorem} {\rm(see \cite{liu3})} \label{t2.1}
Suppose $\gamma \in
\mathcal {P}(2n)$ is a fundamental solution of (\ref{e2.2}) with
$B(t)>0$. There holds
\begin{eqnarray*}
i_{L}(\gamma)\geq 0.
\end{eqnarray*}
\end{theorem}

Suppose the continuous symplectic path $\gamma: [0, 2]\rightarrow
Sp(2n)$ is the fundamental solution of (\ref{e2.2}) with $B(t)$
satisfying $B(t+2)=B(t)$ and $B(1+t)N=NB(1-t)$. This implies
$B(t)N=NB(-t)$. By the unique existence theorem of the differential
equations, we get
\begin{eqnarray*}
\gamma(1+t)=N\gamma(1-t)\gamma(1)^{-1}N\gamma(1),\;
\gamma(2+t)=\gamma(t)\gamma(2).
\end{eqnarray*}
We define the iteration path of $\gamma|_{[0, 1]}$ by
\begin{eqnarray*}
\gamma^{1}(t)&=&\gamma(t),\; t\in [0,1],\\
\gamma^{2}(t)&=&\left\{\begin{array}{ll}&\gamma(t),\; t \in [0, 1], \\
&N\gamma(2-t)\gamma(1)^{-1}N\gamma(1),\; t\in [1,
2],\end{array}\right.\\
\gamma^{3}(t)&=&\left\{\begin{array}{ll}&\gamma(t),\; t \in [0, 1], \\
&N\gamma(2-t)\gamma(1)^{-1}N\gamma(1),\; t\in [1,
2], \\ &\gamma(t-2)\gamma(2),\; t\in [2, 3],\end{array}\right.\\
\gamma^{4}(t)&=&\left\{\begin{array}{ll}&\gamma(t),\; t \in [0, 1], \\
&N\gamma(2-t)\gamma(1)^{-1}N\gamma(1),\; t\in [1, 2], \\
&\gamma(t-2)\gamma(2),\; t\in [2, 3], \\
&N\gamma(4-t)\gamma(1)^{-1}N\gamma(1)\gamma(2),\; t\in [3,
4],\end{array}\right.
\end{eqnarray*}
and in general, for $k\in \mathbb{N}$, we define
\begin{eqnarray*}
\gamma^{2k-1}(t)&=&\left\{\begin{array}{ll}&\gamma(t),\; t \in [0, 1], \\
&N\gamma(2-t)\gamma(1)^{-1}N\gamma(1),\; t\in [1, 2], \\
&\cdots\cdots \\
&N\gamma(2k-2-t)\gamma(1)^{-1}N\gamma(1)\gamma(2)^{2k-5},\; t\in
[2k-3, 2k-2], \\
&\gamma(t-2k+2)\gamma(2)^{2k-4},\; t\in[2k-2, 2k-1],
\end{array}\right.\\
\gamma^{2k}(t)&=&\left\{\begin{array}{ll}&\gamma(t),\; t \in [0, 1], \\
&N\gamma(2-t)\gamma(1)^{-1}N\gamma(1),\; t\in [1, 2], \\
&\cdots\cdots \\
&\gamma(t-2k+2)\gamma(2)^{2k-4},\; t\in[2k-2, 2k-1],\\
&N\gamma(2k-t)\gamma(1)^{-1}N\gamma(1)\gamma(2)^{2k-3},\; t\in
[2k-1, 2k].\end{array}\right.
\end{eqnarray*}
Recall that $(i_{\omega}(\gamma), \nu_{\omega}(\gamma))$ is the
$\omega$-index pair of the symplectic path $\gamma$ introduced in
\cite{lo1}, and $(i_{\omega}^{L_{0}}(\gamma),
\nu_{\omega}^{L_{0}}(\gamma))$ is defined in \cite{liu5}.

\begin{theorem} {\rm(see \cite{liu5})} \label{t2.2}
Suppose
$\omega_{k}=e^{\pi \sqrt{-1}/k}$. For odd $k$ we have
\begin{eqnarray*}
i_{L_{0}}(\gamma^{k})&=&i_{L_{0}}(\gamma^{1})+\sum\limits_{i=1}^{(k-1)/2}
i_{\omega_{k}^{2i}}(\gamma^{2}),\\
\nu_{L_{0}}(\gamma^{k})&=&\nu_{L_{0}}(\gamma^{1})+\sum\limits_{i=1}^{(k-1)/2}
\nu_{\omega_{k}^{2i}}(\gamma^{2}),
\end{eqnarray*}
for even $k$, we have
\begin{eqnarray*}
i_{L_{0}}(\gamma^{k})&=&i_{L_{0}}(\gamma^{1})+i_{\omega_{k}^{k/2}}^{L_{0}}(\gamma^{1})
+\sum\limits_{i=1}^{k/2-1}
i_{\omega_{k}^{2i}}(\gamma^{2}),\\
\nu_{L_{0}}(\gamma^{k})&=&\nu_{L_{0}}(\gamma^{1})+\nu_{\omega_{k}^{k/2}}^{L_{0}}(\gamma^{1})
+\sum\limits_{i=1}^{k/2-1} \nu_{\omega_{k}^{2i}}(\gamma^{2}),
\end{eqnarray*}
where $\omega_{k}^{k/2}=\sqrt{-1}$.
\end{theorem}

\begin{theorem} {\rm(see \cite{liu5})} \label{t2.3}
There hold
\begin{eqnarray*}
i_{1}(\gamma^{2})&=&i_{L_{0}}(\gamma^{1})+i_{L_{1}}(\gamma^{1})+n,\\
\nu_{1}(\gamma^{2})&=&\nu_{L_{0}}(\gamma^{1})+\nu_{L_{1}}(\gamma^{1}),
\end{eqnarray*}
where $L_{1}=\mathbb{R} \oplus \{0\}\in \Lambda(n)$.
\end{theorem}

In the following section, we need the following two iteration
inequalities.

\begin{theorem} {\rm(see \cite{liu2})} \label{t2.4}
For any $\gamma \in
\mathcal {P}(2n)$ and $k\in \mathbb{N}$, there holds
\begin{eqnarray*}
&&i_{L_{0}}(\gamma^{1})+\frac{k-1}{2}
\left(i_{1}(\gamma^{2})+\nu_{1}(\gamma^{2})-n\right)\leq
i_{L_{0}}(\gamma^{k})\\
&&\leq
i_{L_{0}}(\gamma^{1})+\frac{k-1}{2}\left(i_{1}(\gamma^{2})+n\right)
-\frac{1}{2}\nu_{1}(\gamma^{2k})+\frac{1}{2}\nu_{1}(\gamma^{2}),\;
\mbox{if }k\in 2\mathbb{N}-1,\\
&&i_{L_{0}}(\gamma^{1})+i_{\sqrt{-1}}^{L_{0}}(\gamma^{1})
+\left(\frac{k}{2}-1\right)\left(i_{1}(\gamma^{2})+\nu_{1}(\gamma^{2})-n\right)\leq
i_{L_{0}}(\gamma^{k})\leq
i_{L_{0}}(\gamma^{1})+i_{\sqrt{-1}}^{L_{0}}(\gamma^{1})\\
&&+\left(\frac{k}{2}-1\right)\left(i_{1}(\gamma^{2})+n\right)
-\frac{1}{2}\nu_{1}(\gamma^{2k})+\frac{1}{2}\nu_{1}(\gamma^{2})
+\frac{1}{2}\nu_{-1}(\gamma^{2}),\; \mbox{if }k\in 2\mathbb{N}.
\end{eqnarray*}
\end{theorem}

\begin{remark}\label{r2.1}
{\rm From (3.17) of \cite{liu5} and Proposition B of \cite{lo2}, we
have that
\begin{eqnarray*}
&i_{L_{0}}(B)\leq i^{L_{0}}_{\omega}(B)\leq i_{L_{0}}(B)+n,\\
&|i_{L_{0}}(B)-i_{L_{1}}(B)|\leq n,
\end{eqnarray*}
where $L_{1}=\mathbb{R} \oplus \{0\}\in \Lambda(n)$.}
\end{remark}

\section{\large\textbf{Proof of Theorem \ref{t1.1} and \ref{t1.2}}}
\setcounter{equation}{0}
\renewcommand{\theequation}{\arabic{section}.\arabic{equation}}

\qquad In this section, we first consider the following Hamiltonian
systems
\begin{eqnarray}\label{e3.1}
 \left\{\begin{array}{ll}&\dot{z}(t)=J\nabla
H(t,z(t)),\; \forall z
\in \mathbb{R}^{2n},\; \forall t \in [0, jT/2], \\
&z(0)\in L_{0},\; z(jT/2)\in L_{0},
 \end{array}\right.
 \end{eqnarray}
where $j\in \mathbb{N}$. The following result is the first part of
Theorem \ref{t1.1}.

\begin{theorem} \label{t3.1}
Suppose $H(t,z)\in C^{2}(\mathbb{R} \times {\mathbb{R}}^{2n},
\mathbb{R})$ satisfies (H4)-(H7), then for $1\leq j<2\pi/\beta_{0}
T$, (\ref{e3.1}) possesses at least one nontrivial solution $z_j$
whose $L_{0}$-index pair $(i_{L_{0}}(z_j), \nu_{L_{0}}(z_j))$
satisfies
\begin{eqnarray*}
i_{L_{0}}(z_j)\leq 1 \leq i_{L_{0}}(z_j)+\nu_{L_{0}}(z_j).
\end{eqnarray*}
So we get a nonconstant brake solution $(\tilde{z}_j, jT)$ with
brake period $jT$ of the Hamiltonian system (\ref{e1.1}) by Lemma
\ref{l1.1}.
\end{theorem}

In order to prove Theorem \ref{t3.1}, we need the following
arguments. For simplicity, we suppose $T=2$. Let $X:=\left\{z \in
W^{1/2,2}([0, j],{\mathbb{R}}^{2n}) |z=\sum\limits_{l\in
{\mathbb{Z}}}e^{\frac{l\pi}{j}Jt}z_{l},\,z_{l}\in L_{0},
\,\|z\|_{X}<+\infty\right\}$ be the Hilbert space with the inner
product $${(u, v)}_{X}=j(u_{0}, v_{0})+j\sum\limits_{l\in
\mathbb{Z}} |l|(u_{l}, v_{l}), \quad \forall u,\,v\in X.$$ In the
following, we use $\langle \cdot , \cdot \rangle$ and $\| \cdot \|$
to denote the inner product and norm in $X$, respectively. It is
well known that if $r\in [1,+\infty)$ and $z\in L^{r}([0,
j],{\mathbb{R}}^{2n})$ then there exists a constant $c_{r}> 0$ such
that $\| z\|_{L^{r}}\leq c_{r}\| z\|.$

We define the linear operators $A$ and $\hat{B}$ on $X$ by extending
the bilinear form
\begin{eqnarray*}
\langle Au,v\rangle=\int_{0}^{j} (-J \dot{u}, v)dt,\quad\quad
\langle \hat{B}u,v\rangle=\int_{0}^{j} (\hat{B}(t)u, v)dt.
\end{eqnarray*}
 Then $\hat{B}$ is a compact self-adjoint operator (see \cite{lo1}) and
 $A$ is a self-adjoint operator, i.e.,
 $\langle Au,v\rangle=\langle u,A^{*}v\rangle=\langle u,Av\rangle$.

Indeed, by definition
\begin{eqnarray*}
&&\langle Au,v\rangle = \int_{0}^{j} (-J \dot{u}(t),
v(t))dt=(-J u(t), v(t))|^{j}_{0}-\int_{0}^{j} (-J u(t), \dot{v}(t))dt \\
&&=(-J u(t), v(t))|_{0}^{j}+\int_{0}^{j} (u(t), -J \dot{v}(t))dt=(-J
u(t), v(t))|_{0}^{j}+\langle u, A v\rangle.
\end{eqnarray*}
Since $(-Ju(t),v(t))|_{0}^{j}=\omega_0(u(T),v(T))-
\omega_0(u(0),v(0))=0$, so $\langle Au, v\rangle=\langle u,
Av\rangle$, i.e., $A$ is a self-adjoint operator.

We take the spaces
\begin{eqnarray*}
X_{m}&=&\left\{z\in X | z=\sum\limits_{l=-m}^{m}
e^{\frac{l\pi}{j}Jt}z_{l},\,z_{l}\in L_{0}\right\},\\
X^{+}&=&\left\{z\in X | z=\sum\limits_{l>0}
e^{\frac{l\pi}{j}Jt}z_{l},\,z_{l}\in L_{0}\right\},\\
X^{-}&=&\left\{z\in X | z=\sum\limits_{l<0}
e^{\frac{l\pi}{j}Jt}z_{l},\,z_{l}\in L_{0}\right\},\\
X^{0}&=&L_{0},
\end{eqnarray*}
and $X_{m}^{+}=X_{m}\cap X^{+}$, $X_{m}^{-}=X_{m}\cap X^{-}$. We
have $X_{m}=X_{m}^{+}\oplus X^{0}\oplus X_{m}^{-}$. We also know
that
\begin{eqnarray}
\langle Az, z\rangle=\frac{\pi}{j}\|z\|^{2},\quad \forall z\in
X_{m}^{+},\label{e3.2}\\
\langle Az, z\rangle=-\frac{\pi}{j}\|z\|^{2},\quad \forall z\in
X_{m}^{-}.\label{e3.3}
\end{eqnarray}
Equalities (\ref{e3.2}) and (\ref{e3.3}) can be proved by definition
and direct computation. Let $P_{m}:X \rightarrow X_{m}$ be the
corresponding orthogonal projection for $m\in \mathbb{N}$. Then
$\Gamma=\{P_{m};\; m\in \mathbb{N}\}$ is a Galerkin approximation
scheme with respect to $A$ (see \cite{liu1}).

For any Lagrangian subspace $L\in \Lambda(n)$, suppose $P\in
Sp(2n)\cap O(2n)$ such that $L=PL_{0}$. Then we define $X_{L}=PX$
and $X_{L}^{m}=PX_{m}$. Let $P^{m}:X_{L} \rightarrow X_{L}^{m}$.
Then as above, $\bar{\Gamma}=\{P^{m};\; m\in \mathbb{N}\}$ is a
Galerkin approximation scheme with respect to $A$. For $d>0$, we
denote by $M_{d}^{*}(Q)$, $*=+,0,-$, the eigenspaces corresponding
to the eigenvalues $\lambda$ of the linear operator $Q:X_{L}\to
X_{L}$ belonging to $[d,+\infty)$, $(-d,d)$ and $(-\infty,-d]$,
respectively. And denote by $M^{*}(Q)$, $*=+,0,-$, the eigenspaces
corresponding to the eigenvalues $\lambda$ of $Q$ belonging to
$(0,+\infty)$, $\{0\}$ and $(-\infty,0)$, respectively. For any
adjoint operator $Q$, we denote $Q^{\sharp}=(Q|_{ImQ})^{-1}$, and we
also denote $P^{m}QP^{m}=(P^{m}QP^{m})|_{X_{L}^{m}}$. The following
result is the well known Galerkin approximation formulas, it is
proved in \cite{liu1}.

\begin{theorem} \label{t3.2}
For any $B(t)\in
C([0,1],\mathscr{L}_{s}(\mathbb{R}^{2n}))$ with its the $L$-index
pair $(i_{L}(B), \nu_{L}(B))$ and any constant $0<d\leq
\displaystyle\frac{1}{4}\|(A-B)^{\sharp}\|^{-1}$, there exists
$m_{0}>0$ such that for $m\geq m_{0}$, we have
\begin{eqnarray*}
&\dim M_{d}^{+}(P^{m}(A-B)P^{m})&=mn-i_{L}(B)-\nu_{L}(B), \\
&\dim M_{d}^{-}(P^{m}(A-B)P^{m})&=mn+i_{L}(B)+n, \\
&\dim M_{d}^{0}(P^{m}(A-B)P^{m})&=\nu_{L}(B).
\end{eqnarray*}
\end{theorem}

We need to truncate the function $\hat{H}$ at infinite. That is to
replace $\hat{H}$ by a modified function $\hat{H}_{K}$ which grows
at a prescribed rate near $\infty$. The truncated function was
defined by P. Rabinowitz in \cite{ra1}. Let $K> 0$ and select
$\chi\in C^{\infty}(\mathbb{R},\mathbb{R})$ such that $\chi(s)=1$
for $s\leq K$, $\chi(s)=0$ for $s\geq K+1$, and $\chi'(s)< 0$ for
$s\in (K, K+1)$. Set
\begin{eqnarray*}
\hat{H}_{K}(t,z)=\chi(|z|)\hat{H}(t, z)+(1-\chi(|z|))r_{K}|z|^{4},
\end{eqnarray*}
where $r_{K}=\max\left\{\displaystyle\frac{\hat{H}(t, z)}{|z|^{4}} |
K\leq |z|\leq K+1, t\in [0,j]\right\}.$ It is known that
$\hat{H}_{K}$ still satisfies (H4)-(H6) with $\theta$ being replaced
by $\hat{\theta}=\max\{\theta,1/4\}$, and $|\nabla
\hat{H}_{K}(t,z)|\leq (z,\nabla \hat{H}_{K}(t,z))+b$, where $b>0$ is
a constant.

Define a functional $\varphi$ on $X$ by
\begin{eqnarray*}
\varphi (z)&=&\frac{1}{2} \langle Az,z\rangle-\int_{0}^{j}
H_{K}(t,z(t))dt\\
 &=&\frac{1}{2}
\langle Az,z\rangle-\frac{1}{2} \langle
\hat{B}z,z\rangle-\int_{0}^{j} \hat{H}_{K}(t,z(t))dt.
\end{eqnarray*}

Suppose $W$ is a real Banach space, $g \in C^{1}(W, \mathbb{R})$.
$g$ is said satisfying the (PS) condition, if for any sequence
$\{x_{q}\}\subset W$ satisfying $g(x_{q})$ is bounded and
$g'(x_{q})\rightarrow 0$ as $q\rightarrow \infty$, there exists a
convergent subsequence $\{x_{q_{h}}\}$ of $\{x_{q}\}$ (see
\cite{ra1}). Let $\varphi_{m}=\varphi |_{X_{m}}$ be the restriction
of $\varphi$ on $X_{m}$. Similar to Proposition A of \cite{ba}, we
have the following two lemmas.

\begin{lemma} \label{l3.1}
For all $m\in \mathbb{N}$, $\varphi_{m}$ satisfies the (PS)
condition on $X_{m}$.
\end{lemma}

\begin{lemma} \label{l3.2}
$\varphi$ satisfies the (PS)$^{*}$ condition
on $X$ with respect to $\{z_{m}\}$, i.e., for any sequence
$\{z_{m}\} \subset X$ satisfying $z_{m}\in X_{m}$,
$\varphi_{m}(z_{m})$ is bounded and $\|\varphi_{m} '
(z_{m})\|_{(X_{m})'}\rightarrow 0$ in $(X_{m})'$ as $m\rightarrow
+\infty$, where $(X_{m})'$ is the dual space of $X_{m}$, there
exists a convergent subsequence $\{z_{m_{h}}\}$ of $\{z_{m}\}$ in
$X$.
\end{lemma}

In order to prove Theorem \ref{t3.1}, we need the following
definition and the saddle-point theorem.

\begin{defi} {\rm(see \cite{gho})} \label{d3.1}
Let $E$ be a
$C^{2}$-Riemannian manifold and $D$ be a closed subset of $E$. A
family $\phi(\alpha)$ of subsets of $E$ is said to be a homological
family of dimensional $q$ with boundary $D$ if for some nontrivial
class $\alpha \in H_{q}(E,D)$. The family $\phi(\alpha)$ is defined
by
\begin{eqnarray*}
\phi(\alpha)=\{G\subset E: \alpha \mbox{ is in the image of } i_{*}:
H_{q}(G,D)\rightarrow H_{q}(E,D)\},
\end{eqnarray*}
where $i_{*}$ is the homomorphism induced by the immersion $i:
G\rightarrow E$.
\end{defi}

\begin{theorem} {\rm(see \cite{gho})} \label{t3.3}
For above $E$, $D$ and
$\alpha$, let $\phi(\alpha)$ be a homological family of dimension
$q$ with boundary $D$. Suppose that $f\in C^{2}(E,\mathbb{R})$
satisfies the (PS) condition. Define
\begin{eqnarray*}
c=\inf\limits_{G\in \phi(\alpha)} \sup\limits_{x\in G} f(x).
\end{eqnarray*}
Suppose that $\sup\limits_{x\in D} f(x)< c$ and $f'$ is Fredholm on
\begin{eqnarray*}
\mathscr{K}_{c}(f)\equiv \{x\in E: f'(x)=0,f(x)=c\}.
\end{eqnarray*}
Then there exists an $x\in \mathscr{K}_{c}(f)$ such that the Morse
index $m^{-}(x)$ and the nullity $m^{0}(x)$ of the functional $f$ at
$x$ satisfy
\begin{eqnarray*}
m^{-}(x)\leq q \leq m^{-}(x)+m^{0}(x).
\end{eqnarray*}
\end{theorem}

It is clear that a critical point of $\varphi$ is a solution of
(\ref{e3.1}). For a critical point $z=z(t)$, we define the
linearized systems at $z(t)$ by
\begin{eqnarray*}
\dot{y}(t)=JH''(t,z(t))y(t).
 \end{eqnarray*}
Let $B(t)=H''(t,z(t))$. Then the $L_{0}$-index pair of $z$ is
defined by
$(i_{L_{0}}(z),\nu_{L_{0}}(z))=(i_{L_{0}}(B),\nu_{L_{0}}(B))$.

\noindent{\textbf{Proof of Theorem \ref{t3.1}.}}\quad We follow the
ideas of \cite{li1} to prove Theorem \ref{t3.1}. We carry out the
proof in 3 steps.

{\textbf{Step 1}} The critical points of $\varphi_{m}$.

Set $S_{m}=X_{m}^{-}\oplus X^{0}$. Then $\mbox{dim }
S_{m}=mn+\mbox{dim }X^{0}=mn+\mbox{dim }\mbox{ker }A=mn+n$,
$\mbox{dim }X_{m}^{+}=mn$.

In the following, we prove that $\varphi_{m}(z)$ satisfies:

(I) $\varphi_{m}(z)\geq \beta>0$, $\forall z\in Y_{m}=X_{m}^{+} \cap
\partial {B_{\rho}(0)}$,

(II) $\varphi_{m}(z)\leq 0<\beta$, $\forall z\in
\partial Q_{m}$, where $Q_{m}=\{re|r\in
 [0,r_{1}]\}\oplus ({B_{r_{2}}(0)}\cap S_{m}),$ $e\in X^{+}_{m}\cap \partial {B_{1}(0)},$
 $r_{1}>\rho,\,r_{2}> 0$.

 First we prove (I). By (H5), for any $\varepsilon> 0$, there is a
$\delta> 0$ such that $\hat{H}_{K}(t,z)\leq \varepsilon|z|^{2}$ if
$|z|\leq \delta$. Since $\hat{H}_{K}(t,z)|z|^{-4}$ is uniformly
bounded as $|z|\rightarrow +\infty$, there is an $M_{1}=M_{1}(K)$
such that $\hat{H}_{K}(t,z)\leq M_{1}|z|^{4}$ for $|z|\geq \delta$.
Hence
\begin{eqnarray*}
\hat{H}_{K}(t,z)\leq \varepsilon|z|^{2}+M_{1}|z|^{4},\quad \forall
z\in \mathbb{R}^{2n}.
\end{eqnarray*}
For $z \in Y_{m}$, we have
\begin{eqnarray}\label{e3.4}
\int_{0}^{j} \hat{H}_{K}(t,z)dt\leq
\varepsilon\|z\|^{2}_{L^{2}}+M_{1}\|z\|^{4}_{L^{4}}\leq (\varepsilon
c_{2}^{2}+M_{1}c_{4}^{4}\|z\|^{2})\|z\|^{2}.
\end{eqnarray}
By (\ref{e3.2}) and (\ref{e3.4})
\begin{eqnarray*}
\varphi_{m}(z) &=& \frac{1}{2} \langle Az,z \rangle
-\frac{1}{2}\langle \hat{B}z, z\rangle-\int_{0}^{j} \hat{H}_{K}(t,z(t))dt \\
 &\geq&
 \frac{\pi}{2j}\|z\|^{2}-\frac{\beta_0}{2}\|z\|^{2}-
(\varepsilon c_{2}^{2}+M_{1}c_{4}^{4}\|z\|^{2})\|z\|^{2} \\
&=& \frac{\pi}{2j}\rho^{2}-\frac{\beta_0}{2}\rho^{2}-(\varepsilon
c_{2}^{2}+M_{1}c_{4}^{4}\rho^{2})\rho^{2}.
\end{eqnarray*}
Since $1\leq j<\pi/\beta_0$, we can choose constants
$\rho=\rho(K)>0$ and $\beta=\beta(K)>0$, which are sufficiently
small and independent of $m$, such that for $z\in Y_{m}$,
\begin{eqnarray*}
\varphi_{m}(z)\geq \beta>0.
\end{eqnarray*}
Hence (I) holds.

Next prove (II). Let $e\in X_{m}^{+}\cap \partial B_{1}$ and
$z=z^{-}+z^{0}\in S_{m}$. By (\ref{e3.2}) and (\ref{e3.3}), there
holds
\begin{eqnarray}
\varphi_{m}(z+re) &=& \frac{1}{2}\langle Az^{-}, z^{-}
\rangle+\frac{1}{2}r^{2}\langle
Ae, e \rangle-\langle \hat{B}(z+re), z+re\rangle-\int_{0}^{j} \hat{H}_{K}(t,z+re)dt\nonumber \\
&\leq& -\frac{\pi}{2j}\|z^{-}\|^{2}+
\frac{\pi}{2j}r^{2}-\int_{0}^{j}\hat{H}_{K}(t,z+re)dt,\label{e3.5}
\end{eqnarray}
If $r=0$, by (H4), we see that
\begin{eqnarray}\label{e3.6}
\varphi_{m}(z+re)\le  -\frac{\pi}{2j}\|z^{-}\|^{2}\le
0.\end{eqnarray}
 If $r=r_1$, or $\|z\|=r_2$, then from (H6), We have
\begin{eqnarray}\label{e3.7}
\hat{H}_{K}(t,z)\geq b_{1}|z|^{\frac{1}{\hat{\theta}}}-b_{2},
\end{eqnarray}
 where
$b_{1}>0,\,b_{2}$ are two constants independent of $K$ and $m$. Then
by (\ref{e3.7}),
\begin{eqnarray}
\int_{0}^{j} \hat{H}_{K}(t,z+re)dt &\geq& b_{1}\int_{0}^{j}
|z+re|^{\frac{1}{\hat{\theta}}}dt-jb_{2}\nonumber \\ &\geq&
b_{3}\left(\int_{0}^{j}
|z+re|^{2}dt\right)^{\frac{1}{2\hat{\theta}}}-b_{4}\nonumber \\
&\geq&
b_{5}\left(\|z^{0}\|^{\frac{1}{\hat{\theta}}}+r^{\frac{1}{\hat{\theta}}}\right)-b_{4},\label{e3.8}
\end{eqnarray}
where $b_{3},\,b_{4}$ are constants and $b_{5}>0$ independent of $K$
and $m$. Thus by (\ref{e3.8}), (\ref{e3.5}) is
\begin{eqnarray}\label{e3.9}
\varphi_{m}(z+re)\leq -\frac{\pi}{2j}\|z^{-}\|^{2}+
\frac{\pi}{2j}r^{2}-
b_{5}\left(\|z^{0}\|^{\frac{1}{\hat{\theta}}}+r^{\frac{1}{\hat{\theta}}}\right)+b_{4},
\end{eqnarray}
Thus we can choose large enough $r_{1}$ and $r_{2}$ independent of
$K$ and $m$ such that
\begin{eqnarray*}
\varphi_{m}(z+re)\leq 0, \quad \mbox{ on } \partial Q_{m}.
\end{eqnarray*}
Then (II) holds.

Because $Q_{m}$ is deformation retract of
 $X_{m}$, then $H_{q}(Q_{m},\partial Q_{m})\cong H_{q}(X_{m},\partial
 Q_{m})$, where $q=\mbox{dim }S_{m}+1=mn+n+1=\mbox{dim }Q_{m}$, and $\partial
 Q_{m}$ is the boundary of $Q_{m}$ in $S_{m}\oplus \{\mathbb{R}e\}$.
 But $H_{q}(Q_{m},\partial Q_{m})\cong H_{q-1}(S^{q-1})\cong
 \mathbb{R}$. Denote by $i:Q_{m}\rightarrow X_{m}$ the inclusion
 map. Let $\alpha=[Q_{m}]\in H_{q}(Q_{m},D)$ be a generator. Then
 $i_{*}\alpha$ is nontrivial in $H_{q}(X_{m},\partial Q_{m})$, and
 $\phi (i_{*}\alpha)$ defined by Definition \ref{d3.1} is a homological
 family of dimension $q$ with boundary $D:=\partial Q_{m}$ and $Q_{m}\in
 \phi(i_{*}\alpha)$. $\partial Q_{m}$ and $Y_{m}$ are homologically link (see \cite{ch}).
 By Lemma \ref{l3.1}, $\varphi_{m}$ satisfies the (PS) condition.
 Define $c_{m}=\inf\limits_{G\in \phi(i_{*}\alpha)} \sup\limits_{z\in G}
 \varphi_{m}(z)$. We have
\begin{eqnarray}\label{e3.10}
\sup\limits_{z\in \partial Q_{m}}
 \varphi_{m}(z)\leq 0< \beta \leq c_{m} \leq \sup\limits_{z\in Q_{m}}
 \varphi_{m}(z) \leq \frac{\pi}{2j} r^{2}_{1}.
 \end{eqnarray}
 Since $X_{m}$ is finite dimensional, $\varphi_{m}'$ is Fredholm.
 By Theorem \ref{t3.3}, $\varphi_{m}$ has a critical point $z_{j}^{m}$ with
 critical value $c_{m}$, and the Morse index $m^{-}(z_{j}^{m})$ and
 nullity $m^{0}(z_{j}^{m})$ of $z_{j}^{m}$ satisfy
\begin{eqnarray}\label{e3.11}
m^{-}(z_{j}^{m})\leq mn+n+1\leq
 m^{-}(z_{j}^{m})+m^{0}(z_{j}^{m}).
 \end{eqnarray}

 Since $\{c_{m}\}$ is bounded, passing to a subsequence, suppose
 $c_{m}\rightarrow c\in [\beta,\frac{\pi}{2j} r_{1}^{2}]$. By the
 (PS)$^{*}$ condition of Lemma \ref{l3.2}, passing to a subsequence, there exists an
 $z_{j}\in X$ such that $$z_{j}^{m}\rightarrow
 z_{j},\,\varphi(z_{j})=c,\,\varphi'(z_{j})=0.$$

 {\textbf{Step 2}} The solution of (\ref{e3.1}).

 Because the critical value $c$ has an upper bound
 $\frac{\pi}{2j}r_{1}^{2}$ independent of $K$, then
 \begin{eqnarray}
 \frac{\pi}{2j}r_{1}^{2}&\geq& c=
\varphi(z_{j})-\frac{1}{2}\langle \varphi'(z_{j}), z_{j}\rangle\nonumber \\
&\geq& \left(\frac{1}{2}-\hat{\theta}\right)\int_{0}^{j}
(z_{j},\nabla \hat{H}_{K}(t,z_{j}))dt.\label{e3.12}
 \end{eqnarray}
 Then by (\ref{e3.12}), $\int_{0}^{j}
(z_{j},\nabla \hat{H}_{K}(t,z_{j}))dt$ has an upper bound
independent of $K$,
\begin{eqnarray}\label{e3.13}
\int_{0}^{j} (z_{j},\nabla \hat{H}_{K}(t,z_{j}))dt\leq \bar{M},
\mbox{ for constant $\bar{M}$ independent of $K$}.
\end{eqnarray}
By $\hat{H}_{K}(t,z_{j})\leq \hat{\theta}(z_{j},\nabla
\hat{H}_{K}(t,z_{j}))$, then by (\ref{e3.13}),
\begin{eqnarray}\label{e3.14}
\int_{0}^{j} \hat{H}_{K}(t,z_{j})dt\leq \hat{\theta}\bar{M}.
\end{eqnarray}
Thus by (\ref{e3.7}) and (\ref{e3.14}),
\begin{eqnarray*}
\int_{0}^{j}
\left(b_{1}|z_{j}|^{\frac{1}{\hat{\theta}}}-b_{2}\right)dt\leq
\int_{0}^{j} \hat{H}_{K}(t,z_{j})dt\leq \hat{\theta}\bar{M},
\end{eqnarray*}
 i.e.,
\begin{eqnarray}\label{e3.15}
\hat{\theta}\bar{M}\geq b_{1}\int_{0}^{j}
|z_{j}|^{\frac{1}{\hat{\theta}}}dt-b_{2}j\geq
b'_{1}\left(\int_{0}^{j}
|z_{j}|^{2}dt\right)^{\frac{1}{2\hat{\theta}}}-b_{2}j.
\end{eqnarray}
Thus by (\ref{e3.15}),
$\|z_{j}\|_{L^{2}}^{\frac{1}{\hat{\theta}}}\leq M_{2}$, where
$M_{2}$ is independent of $K$, i.e.,
\begin{eqnarray}\label{e3.16}
\|z_{j}\|_{L^{2}}\leq M_{3}, \mbox{ where $M_{3}$ is independent of
$K$}.
\end{eqnarray}
  Since
\begin{eqnarray}\label{e3.17}
\|z_{j}\|_{L^{1}}\leq C\|z_{j}\|_{L^{2}}\leq M_{3}', \mbox{ where
$C>0$ is independent of $K$}.
\end{eqnarray}
Thus by (\ref{e3.16}) and (\ref{e3.17}), $\|z_{j}\|_{L^{1}}$ has an
upper bound independent of $K$. We use Young's inequality. For any
$w\in W^{1,2}([0,j], \mathbb{R}^{2n}),\,w(\tau)-w(t)=\int_{t}^{\tau}
\dot{w}(s)ds$. Integrating with respect to $t$ shows that
\begin{eqnarray*}
jw(\tau)-\int_{0}^{j} w(t)dt=\int_{0}^{j} \int_{t}^{\tau}
\dot{w}(s)dsdt,
\end{eqnarray*}
  i.e.,
\begin{eqnarray*}
|jw(\tau)| &=& \left|\int_{0}^{j} w(t)dt+\int_{0}^{j}
\int_{t}^{\tau} \dot{w}(s)dsdt\right|\\ &\leq& \int_{0}^{j}
|w(t)|dt+\int_{0}^{j} \int_{t}^{\tau} |\dot{w}(s)|dsdt\\  &\leq&
\|w\|_{L^{1}}+j\int_{0}^{j} |\dot{w}(s)|ds\\ &=&
\|w\|_{L^{1}}+j\|\dot{w}\|_{L^{1}}, \end{eqnarray*} i.e.,
\begin{eqnarray*}
|w(\tau)|\leq \frac{\|w\|_{L^{1}}}{j}+\|\dot{w}\|_{L^{1}},
\end{eqnarray*}
 i.e.,
\begin{eqnarray}\label{e3.18}
\|w\|_{L^{\infty}}\leq
\displaystyle\frac{\|w\|_{L^{1}}}{j}+\|\dot{w}\|_{L^{1}}.
\end{eqnarray}
Therefore
\begin{eqnarray}
\|\dot{z}_{j}\|_{L^{1}}&=& \|J\hat{B}(t)z_{j}+J\nabla
\hat{H}_{K}(t,z_{j})\|_{L^{1}}\nonumber \\ &\leq&
\beta_0\|z_{j}\|_{L^{1}}+\int_{0}^{j} |\nabla \hat{H}_{K}(t,z_{j})|dt\nonumber \\
&\leq& \beta_0 M'_{3}+\int_{0}^{j} (z_{j},\nabla
\hat{H}_{K}(t,z_{j}))dt+bj\nonumber
\\ &\leq& \beta_0 M'_{3}+\bar{M}j+bj \leq M_{4}.\label{e3.19}
\end{eqnarray}
Thus $\|\dot{z}_{j}\|_{L^{1}}$ has an upper bound independent of
$K$. Then from (\ref{e3.17}), (\ref{e3.18}) and (\ref{e3.19}),
$\|z_{j}\|_{L^{\infty}}\leq K_0$, where $K_{0}$ is independent of
$K$. We choose $K>K_{0}$, therefore
$\hat{H}_{K}(t,z_{j})=\hat{H}(t,z_{j})$. Consequently, $z_{j}$ is a
nontrivial solution of (\ref{e3.1}). Then by Lemma \ref{l1.1}, we
get a nonconstant brake solution $\tilde{z}_j$ of the Hamiltonian
system (\ref{e1.1}).

{\textbf{Step 3}}\quad Let $B(t)=H_{K}''(t,z_{j}(t))$,
 $d=\displaystyle\frac{1}{4}\|(A-B)^{\sharp}\|^{-1}.$ Since
 \begin{eqnarray*}
 \|\varphi''(x)-(A-B)\|\rightarrow 0 \quad \mbox{as } \|x-z_{j}\|\rightarrow 0,
 \end{eqnarray*}
there exists a $r_{3}>0$ such that
\begin{eqnarray*}
 \|\varphi''(x)-(A-B)\|<\frac{1}{4} d,\quad \forall x\in
 V_{r_{3}}(z_{j})=\{x\in X |\quad \|x-z_{j}\|\leq r_{3}\}.
 \end{eqnarray*}
 Then for $m$ large enough, there holds
 \begin{eqnarray}\label{e3.20}
 \|\varphi''_{m}(x)-P_{m}(A-B)P_{m}\|<\frac{1}{2} d,\quad \forall x\in
 V_{r_{3}}(z_{j})\cap X_{m}.
 \end{eqnarray}
 For $x\in V_{r_{3}}(z_{j})\cap X_{m}$, $\forall u\in
 M^{-}_{d}(P_{m}(A-B)P_{m})\setminus \{0\}$, from (\ref{e3.20}) we have
\begin{eqnarray}
 \langle \varphi''_{m}(x)u,u \rangle &\leq& \langle
 P_{m}(A-B)P_{m}u,u
 \rangle+\|\varphi''_{m}(x)-P_{m}(A-B)P_{m}\|\cdot \|u\|^{2}\nonumber \\
 & \leq & -\frac{1}{2} d \|u\|^{2}<0.\label{e3.21}
 \end{eqnarray}
 Thus by (\ref{e3.21}),
\begin{eqnarray}\label{e3.22}
 \dim M^{-}(\varphi''_{m}(x)) \geq \dim
 M^{-}_{d}(P_{m}(A-B)P_{m}),\quad \forall x\in V_{r_{3}}(z_{j})\cap X_{m}.
 \end{eqnarray}
 Similarly, we have
\begin{eqnarray}\label{e3.23}
\dim M^{+}(\varphi''_{m}(x)) \geq \dim
 M^{+}_{d}(P_{m}(A-B)P_{m}),\quad \forall x\in V_{r_{3}}(z_{j})\cap X_{m}.
 \end{eqnarray}
 By (\ref{e3.11}), (\ref{e3.22}), (\ref{e3.23}) and Theorem \ref{t3.2}, for large $m$ we have
 \begin{eqnarray}
 mn+n+1 & \geq & m^{-}(z_{j}^{m})\nonumber \\ & \geq & \dim
 M^{-}_{d}(P_{m}(A-B)P_{m})\nonumber \\ &=&
 mn+i_{L_{0}}(B)+n.\label{e3.24}
 \end{eqnarray}
 We also have
 \begin{eqnarray}
 mn+n+1 & \leq & m^{-}(z_{j}^{m})+m^{0}(z_{j}^{m})\nonumber \\ & \leq & \dim
 M^{-}_{d}(P_{m}(A-B)P_{m})\oplus \dim M^{0}_{d}(P_{m}(A-B)P_{m})\nonumber \\
 &=& mn+i_{L_{0}}(B)+n+\nu_{L_{0}}(B).\label{e3.25}
 \end{eqnarray}
Combining (\ref{e3.24}) and (\ref{e3.25}), we have
\begin{eqnarray*}
 i_{L_{0}}(z_{j})\leq 1 \leq
 i_{L_{0}}(z_{j})+\nu_{L_{0}}(z_{j}).
 \end{eqnarray*}

The proof of Theorem \ref{t3.1} is complete.$\hfill\Box$

It is the time to give the proof of Theorem \ref{t1.1} and
\ref{t1.2}.

\noindent{\textbf{Proof of Theorem \ref{t1.1}.}}\quad For $1\leq k<
\pi/\beta_0$, by Theorem \ref{t3.1}, we obtain that there is a
nontrivial solution $(z_{k}, k)$ of the Hamiltonian systems
(\ref{e3.1}) and its $L_{0}$-index pair satisfies
\begin{eqnarray}\label{e3.26}
 i_{L_{0}}(z_{k}, k)\leq 1 \leq
i_{L_{0}}(z_{k}, k)+\nu_{L_{0}}(z_{k}, k).
\end{eqnarray}
Then by Lemma \ref{l1.1}, $(\tilde{z}_{k}, 2k)$ is a nonconstant
brake solution of (\ref{e1.1}).

For $k\in 2\mathbb{N}-1$, we suppose that $(\tilde{z}_{1}, 2)$ and
$(\tilde{z}_{k}, 2k)$ are not distinct. By (\ref{e3.26}), Theorem
\ref{t2.3} and Theorem \ref{t2.4}, we have
\begin{eqnarray}
&&1 \geq i_{L_{0}}(z_{k}, k) \geq i_{L_{0}}(z_{1}, 1)+\frac{k-1}{2}
\left(i_{1}(\tilde{z}_{1}, 2)+\nu_{1}(\tilde{z}_{1},
2)-n\right)\nonumber\\
&\geq& i_{L_{0}}(z_{1}, 1)+\frac{k-1}{2} \left(i_{L_{0}}(z_{1},
1)+i_{L_{1}}(z_{1}, 1)+n+\nu_{L_{0}}(z_{1},
1)+\nu_{L_{1}}(z_{1}, 1)-n\right)\nonumber\\
&=&i_{L_{0}}(z_{1}, 1)+\frac{k-1}{2}\left(i_{L_{0}}(z_{1},
1)+i_{L_{1}}(z_{1}, 1)+\nu_{L_{0}}(z_{1}, 1)+\nu_{L_{1}}(z_{1},
1)\right),\label{e3.27}
\end{eqnarray}
where $L_{1}=\mathbb{R}^{n} \oplus \{0\} \in \Lambda (n)$. By (H3),
(H7) and Theorem \ref{t2.1}, we have $i_{L_{1}}(z_{1}, 1)\geq 0$. We
also know that $\nu_{L_{1}}(z_{1}, 1)\geq 0$ and $i_{L_{0}}(z_{1},
1)+\nu_{L_{0}}(z_{1}, 1)\geq 1$. Then (\ref{e3.27}) is
\begin{eqnarray}\label{e3.28}
1 \geq i_{L_{0}}(z_{1}, 1)+\frac{k-1}{2}.
\end{eqnarray}
By $0\leq i_{L_{0}}(z_{1}, 1)\leq 1$, from (\ref{e3.28}) we have
$\frac{k-1}{2} \leq 1$, i.e., $k\leq 3$. It is contradict to $k\geq
5$. Similarly, we have that for each $k\in 2\mathbb{N}-1$, $k\geq 5$
and $kj< \frac{\pi}{\beta_0}$, $1\leq j< \frac{\pi}{\beta_0}$,
$(\tilde{z}_{j}, 2j)$ and $(\tilde{z}_{kj}, 2kj)$ are distinct brake
solutions of (\ref{e1.1}). Furthermore, $(\tilde{z}_{1}, 2)$,
$(\tilde{z}_{k}, 2k)$, $(\tilde{z}_{k^{2}}, 2k^{2})$,
$(\tilde{z}_{k^{3}}, 2k^{3})$, $\cdots$, $(\tilde{z}_{k^{p}},
2k^{p})$ are pairwise distinct brake solutions of (\ref{e1.1}),
where $k\in 2\mathbb{N}-1$, $k\geq 5$ and $1\leq k^{p}<
\frac{\pi}{\beta_0}$ with $p\in \mathbb{N}$.

For $k\in 2\mathbb{N}$, as above, we suppose that $(\tilde{z}_{1},
2)$ and $(\tilde{z}_{k}, 2k)$ are not distinct. By (\ref{e3.26}),
Theorem \ref{t2.3} and Theorem \ref{t2.4}, we have
\begin{eqnarray}
&&1 \geq i_{L_{0}}(z_{k}, k) \geq i_{L_{0}}(z_{1},
1)+i^{L_{0}}_{\sqrt{-1}}(z_{1},
1)+\left(\frac{k}{2}-1\right)\left(i_{1}(\tilde{z}_{1},
2)+\nu_{1}(\tilde{z}_{1},
2)-n\right)\nonumber\\
&\geq& i_{L_{0}}(z_{1}, 1)+i^{L_{0}}_{\sqrt{-1}}(z_{1},
1)+\left(\frac{k}{2}-1\right)(i_{L_{0}}(z_{1},
1)+i_{L_{1}}(z_{1}, 1)+n\nonumber\\
&&+\nu_{L_{0}}(z_{1},
1)+\nu_{L_{1}}(z_{1}, 1)-n)\nonumber\\
&=& i_{L_{0}}(z_{1}, 1)+i^{L_{0}}_{\sqrt{-1}}(z_{1},
1)+\left(\frac{k}{2}-1\right) (i_{L_{0}}(z_{1}, 1)+i_{L_{1}}(z_{1},
1)+\nu_{L_{0}}(z_{1},
1)\nonumber\\
&&+\nu_{L_{1}}(z_{1}, 1)).\label{e3.29}
\end{eqnarray}
Similarly, we also know that $i_{L_{1}}(z_{1}, 1)\geq 0$,
$\nu_{L_{1}}(z_{1}, 1)\geq 0$, $i_{L_{0}}(z_{1},
1)+\nu_{L_{0}}(z_{1}, 1)\geq 1$. By Remark \ref{r2.1}, we have
$i^{L_{0}}_{\sqrt{-1}}(z_{1}, 1)\geq i_{L_{0}}(z_{1}, 1)\geq 0$.
Then (\ref{e3.29}) is
\begin{eqnarray}\label{e3.30}
1 \geq i_{L_{0}}(z_{1}, 1)+\left(\frac{k}{2}-1\right).
\end{eqnarray}
By $0\leq i_{L_{0}}(z_{1}, 1)\leq 1$,  from (\ref{e3.30}) we have
$\frac{k}{2}-1 \leq 1$, i.e., $k\leq 4$. It  contradicts  to $k\geq
5$. Similarly we have that for each $k\in 2\mathbb{N}$, $k\geq 6$
and $kj< \frac{\pi}{\beta_0}$, $1\leq j< \frac{\pi}{\beta_0}$,
$(\tilde{z}_{j}, 2j)$ and $(\tilde{z}_{kj}, 2kj)$ are distinct brake
solutions of (\ref{e1.1}). Furthermore, $(\tilde{z}_{1}, 2)$,
$(\tilde{z}_{k}, 2k)$, $(\tilde{z}_{k^{2}}, 2k^{2})$,
$(\tilde{z}_{k^{3}}, 2k^{3})$, $\cdots$, $(\tilde{z}_{k^{p}},
2k^{p})$ are pairwise distinct brake solutions of (\ref{e1.1}),
where $k\in 2\mathbb{N}$, $k\geq 6$ and $1\leq k^{p}<
\frac{\pi}{\beta_0}$ with $p\in \mathbb{N}$.

In all, for any integer $1\leq j<\frac{\pi}{\beta_0}$,
$\tilde{z}_{j}$ and $\tilde{z}_{kj}$ are distinct brake solutions of
(\ref{e1.1}) for $k\geq 5$ and $kj< \frac{\pi}{\beta_0}$.
Furthermore, $\{\tilde{z}_{k^{p}}| p\in \mathbb{N}\}$ is a pairwise
distinct brake solution sequence of (\ref{e1.1}) for $k\geq 5$ and
$1\leq k^{p}< \frac{\pi}{\beta_0}$. The proof of Theorem \ref{t1.1}
is complete. $\hfill\Box$

We note that Theorem \ref{t1.2} is a direct consequence of Theorem
\ref{t1.1}.


\begin{thebibliography}{EE}
\baselineskip=18pt



\bibitem {ABL1} Ambrosetti, A., Benci, V. and Long, Y., A note on the existence of
multiple brake orbits, {\em Nonlinear Anal. TMA}, {\bfseries 21},
1993, 643-649.

\bibitem{ba} Bahri, A. and Berestycki, H., Forced vibrations of
 superquadratic Hamiltonian systems, {\em Acta Math.}, {\bfseries 152}, 1984, 143--197.

\bibitem{Be1} Benci, V., Closed geodesics for the Jacobi metric and periodic
solutions of prescribed energy of natural Hamiltonian systems, {\em
Ann. I. H. P. Analyse Nonl.}, {\bfseries 1}, 1984, 401-412.

\bibitem{BG} Benci, V. and Giannoni, F., A new proof of the existence of a brake
orbit, In ``Advanced Topics in the Theory of Dynamical Systems".
{\em Notes Rep. Math. Sci. Eng.}, {\bfseries 6}, 1989, 37-49.

\bibitem{Bol} Bolotin, S., Libration motions of natural dynamical
systems, {\em Vestnik Moskov Univ. Ser. I. Mat. Mekh.}, {\bfseries
6}, 1978, 72-77 (in Russian).

\bibitem{BolZ} Bolotin, S. and Kozlov, V. V., Librations with many degrees
of freedom, {\em J. Appl. Math. Mech.}, {\bfseries 42}, 1978,
245-250 (in Russian).

\bibitem{ch} Chang, K., Infinite Dimensional Morse Theory and
Multiple Solution Problems, Birkh\"{a}user Verlag, Basel, Boston,
Berlin, 1993.

\bibitem{ek}  Ekeland, I., Convexity Method in
Hamiltonian Mechanics, Berlin: Springer-Verlag, 1990.

\bibitem{ekho}  Ekeland, I. and Hofer, H., Subharmonics of convex Hamiltonian
systems, {\em Comm. Pure Appl. Math.}, {\bfseries 40}, 1987, 1--37.

\bibitem{gho} Ghoussoub, N., Location, multiplicity and Morse
indices of minimax critical points, {\em J. Reine Angew Math.},
{\bfseries 417}, 1991, 27--76.

\bibitem {GZ1} Gluck, H. and Ziller, W., Existence of periodic
solutions of conservtive systems, {\em Seminar on Minimal
Submanifolds}, Princeton University Press 1983, 65-98.

\bibitem {Gro} Groesen, E. W. C. van, Analytical mini-max methods for Hamiltonian
brake orbits of prescribed energy, {\em J. Math. Anal. Appl.},
{\bfseries 132}, 1988, 1-12.

\bibitem {Ha1} Hayashi, K.,
Periodic solution of classical Hamiltonian systems, {\em Tokyo J.
Math.}, {\bfseries 6}, 1983, 473-486.

\bibitem{li1} Li, C. and Liu, C., Nontrivial solutions of
superquadratic Hamiltonian systems with Lagrangian boundary
conditions and the $L$-index theory, {\em Chinese Ann. Math. Ser.
B}, {\bfseries 29(6)}, 2008, 597--610.

\bibitem{liu1} Liu, C., Asymptotically linear Hamiltonian systems
with Lagrangian boundary conditions, {\em Pacific J. Math.},
{\bfseries 232(1)}, 2007, 233--255.

\bibitem{liu3} Liu, C., Maslov-type index theory for symplectic
paths with Lagrangian boundary conditions, {\em Adv. Non. Stu.},
{\bfseries 7}, 2007, 131--161.

\bibitem{liu2} Liu, C. Minimal period estimates for brake orbits of
nonlinear symmetric Hamiltonian systems, preprint.

\bibitem{liu4} Liu, C., Subharmonic solutions of Hamiltonian
systems, {\em Nonlinear Anal. TMA}, {\bfseries 42}, 2000, 185--198.

\bibitem{liu5} Liu, C. and Zhang, D., An iteration theory of
Maslov-type index associated with a Lagrangian subspace for
symplectic paths and multiplicity of brake orbits in bounded convex
symmetric domains, preprint.

\bibitem{lo1} Long, Y., Index Theory for Symplectic Paths with
Applications, Birkh\"{a}user Verlag, Basel, Boston, Berlin, 2002.

\bibitem{lo4}  Long, Y., Multiple solutions of perturbed superquadratic second order Hamiltonian
systems, {\em Trans. Amer. Math. Soc.}, {\bfseries 311}, 1989,
749--780.

\bibitem{lo6}  Long, Y., Periodic solutions of perturbed superquadratic Hamiltonian systems,
 {\em Ann. Scuola Norm. Sup. Pisa.}, Series 4, {\bfseries 17}, 1990, 35--77.

\bibitem{lo2} Long, Y., Zhang, D. and Zhu, C., Multiple brake
orbits in bounded convex symmetric
 domains, {\em Adv. in Math.}, {\bfseries 203}, 2006, 568--635.

\bibitem{mc} McDuff, D. and Salamon, D., Introduction to
Symplectic Topology, Clarendon Press, Oxford, 1998.

\bibitem{ra1} Rabinowitz, P. H., Minimax methods in critical
point theory with applications to differential equations, CBMS
Regional Conf Ser in Math, 65, Ams, RI, 1986.

\bibitem{Ra1} Rabinowitz, P. H., On the existence of periodic solutions for a
class of symmetric Hamiltonian systems, {\em Nonlinear Anal. TMA},
{\bfseries 11}, 1987, 599-611.

\bibitem{ra2} Rabinowitz, P. H., On subharmonic solutions of
Hamiltonian
 systems, {\em Comm. Pure Appl. Math.}, {\bfseries 33}, 1980, 609--633.

\bibitem{ra4}  Rabinowitz, P. H., Periodic
solution of Hamiltonian systems, {\em Comm. Pure Appl. Math.},
{\bfseries 31}, 1978, 157--184.

\bibitem{se} Seifert, H., Periodische Bewegungen mechanischer
Systeme, {\em Math. Z.}, {\bfseries 51}, 1948, 197--216.

\bibitem{si}  Silva, E. A. B., Subharmonic solutions for subquadratic
Hamiltonian systems, {\em J. Diff. Eq.}, {\bfseries 115}, 1995,
120--145.

\bibitem {Sz} Szulkin, A., An index theory and existence of multiple brake orbits
for star-shaped Hamiltonian systems, {\em Math. Ann.}, {\bfseries
283}, 1989, 241-255.

\bibitem{an} An, T., Subharmonic solutions of Hamiltonian
systems and the Maslov-type index theory, {\em J. Math. Anal.
Appl.}, {\bfseries 331}, 2007, 701--711.




\end{thebibliography}
\end{document}